\documentclass[12pt,a4paper]{article}
\usepackage{amsmath}
\usepackage{amssymb}
\usepackage{amsthm}
\usepackage{amsmath}
\usepackage{amssymb}
\usepackage{bbm}
\usepackage{url}
\usepackage{amscd}
\usepackage{mathrsfs}
\usepackage{graphicx}
\usepackage{xcolor}
\usepackage{float}
\usepackage[margin=0.7in]{geometry}
\usepackage{cite}

\newtheorem{definition}{Definition}

\usepackage{setspace}
\usepackage{graphicx}
\usepackage{amssymb}
\usepackage{amsmath}
\usepackage{amsthm}
\usepackage{tikz-cd}
\usepackage{xcolor}
\usepackage{hyperref}
\usepackage{multirow}
\usepackage{enumitem}
\usepackage{listings}
\usepackage{pgfplots}
\usepackage{wrapfig}
\usepackage{lipsum}
\usepackage[utf8x]{inputenc}
\usepackage{float}
\usepackage [english]{babel}
\usepackage [autostyle, english = american]{csquotes}

\newtheorem{Thm}{Theorem}
\newtheorem{Lem}[Thm]{Lemma}

\def\F{\mathbb{F}}

\def\P{\mathbb{P}}

\def\a{\alpha}

\def\z{\mathbb{P}^{1}(\overline{\mathbb{F}}_{q})}
\def\y{\overline{\mathbb{F}}_{q}}
\def\b{\beta}

\def\f{\mathbb{F}_{q^{2}}}
\def\u{\mu_{q+1}}
\def\g{\mathbb{F}_{q}^{\ast}}

\def\pp{\mathbb{P}^{1}(\mathbb{F}_{q})}
\def\h{\mathbb{F}_{q}}

\title{Some permutation pentanomials over finite fields of even characteristic\thanks{This work was supported by the Research Grants Council (RGC) of the Hong Kong (No. 16307524).}}
\author{Farhana Kousar\thanks{F. Kousar is at the Department of Mathematics, Hong Kong University of Science and Technology, Clear Water Bay, Kowloon, Hong Kong, and also at the department of mathematics, Sardar Bahadur Khan Women's University, Quetta, Pakistan. Email: fkousar@connect.ust.hk and farhanakausar66@gmail.com} and Maosheng Xiong\thanks{M. Xiong is at the Department of Mathematics, Hong Kong University of Science and Technology, Clear Water Bay, Kowloon, Hong Kong. Email: mamsxiong@ust.hk}}
\date{}

\begin{document}
	
	\maketitle
	
	\noindent{{\bf Abstract:}
		In a recent paper \cite{zhang2024more} Zhang et al. constructed 17 families of permutation pentanomials of the form $x^t+x^{r_1(q-1)+t}+x^{r_2(q-1)+t}+x^{r_3(q-1)+t}+x^{r_4(q-1)+t}$ over $\f$ where $q=2^m$. In this paper for 14 of these 17 families we provide a simple explanation as to why they are permutations. We also extend these 14 families into three general classes of permutation pentanomials over $\f$.}
	
	
	\section{Introduction}
Let $p$ be a prime, $m$ a positive integer, $q=p^m$, and $\h$ a finite field of $q$ elements. A polynomial $f(x)\in\h[x]$ is called a permutation polynomial (\textbf{PP} for short) over $\h$ if the associated function $a\mapsto f(a)$ from $\h$ to itself is a permutation. Permutation polynomials have been studied extensively in the past thirty years for their connection to many branches of mathematics and their wide applications in areas like cryptography \cite{laigle2007permutation,sun2005interleavers}, coding theory \cite{muller1986cryptanalysis, rivest1978method} and combinatorial design theory \cite{ding2006family}.

Permutation polynomials of the type $x^{r}H(x^{q-1})$ over $\f$, where $r$ is a positive integer, are of significant interest due to their simple algebraic feature and their connection to Niho exponents \cite{LZ}. This particular type of PPs with two, three, or four terms (known in the literature as "binomials," "trinomials," or "quadrinomials") has been discovered in numerous papers (see, for instance, \cite{ayad2014permutation,BARTOLI2021101781,deng2019more,ding2023determination,gupta2024class,gupta2016some, hou2013class,hou2014determination,hou2015determination,hou2015survey,li2017new,li2018new,li2021several,masuda2009permutation,ozbudak2023classification,zheng2022more} and references therein). On the other hand, partly due to technical complexities, there are only a few paper in the literature on permutation pentanomials (PPs with five terms) of the type $x^{r}H(x^{q-1})$ over $\f$ \cite{li2018new,li2021several,xu2018some,zhang2023more,zhang2024more}.

In a most recent paper \cite{zhang2024more}, the authors studied permutation pentanomials over $\mathbb{F}_{2^{2m}}$ of the form
\begin{eqnarray} \label{1:fx} f(x)&=&x^t+x^{r_1(q-1)+t}+x^{r_2(q-1)+t}+x^{r_3(q-1)+t}+x^{r_4(q-1)+t}
\end{eqnarray}
such that
\begin{eqnarray} \label{1:gcd}
	\gcd\left(x^{r_4}+x^{r_3}+x^{r_2}+x^{r_1}+1, x^{t}+x^{t-r_1}+x^{t-r_2}+x^{t-r_3}+x^{t-r_4}\right)=1.
\end{eqnarray}
The polynomial $f(x)$ in (\ref{1:fx}) can be written as $f(x)=x^tH\left(x^{q-1}\right)$ where $H(x)=1+x^{r_1}+x^{r_2}+x^{r_3}+x^{r_4}$. Searching through numerical data for $4 \le t<100, 1 \le r_i \le t, \forall 1 \le i \le 4$, they found and proved 17 families of such permutation pentanomials $f(x)$ satisfying (\ref{1:fx}) and (\ref{1:gcd}) which we reproduce in Table \ref{1:tab} below (see also \cite[Table 2]{zhang2024more}):
\begin{table}[H]
	\centering
	\resizebox{0.9\textwidth}{!}{%
		\begin{tabular}{ccccccc}
			\hline
			No& $f(x)$&$ \mathrm{PP}$ Conditions & Class&$Q_{1}$  & $Q_{2}$ & $Q_{1}+Q_{2}+1$    \\
			\hline
			1 & $x^{8q+1}+x^{7q+2}+x^{5q+4}+x^{3q+6}+x^9$& $m$ is odd& -- & -- & --  & --   \\
			2* & $x^{10q+1}+x^{9q+2}+x^{3q+8}+x^{q+10}+x^{11}$&$ m \not \equiv 0 \pmod{10}$& $f_{A}$ & $2^{3}$ & $2$  & 11   \\
			3 & $x^{10q+1}+x^{8q+3}+x^{6q+5}+x^{4q+7}+x^{11}$&$\gcd(m,5)=1$ & -- & -- & --  & --  \\
			4* & $x^{12q+1}+x^{9q+4}+x^{8q+5}+x^{5q+8}+x^{13}$& $m \not \equiv 0 \pmod{6}$ & $f_B$& $2^{3}$ & $2^{2}$  & $13$\\
			5*  & $x^{18q+1}+x^{17q+2}+x^{3q+16}+x^{2q+17}+x^{19}$& $m \not \equiv 0 \pmod{18}$&$f_{B}$  &$2$ & $2^{4}$ & $19$  \\
			6*  & $x^{21q}+x^{16q+5}+x^{4q+17}+x^{q+20}+x^{21}$& $m \not \equiv 3 \pmod{6}$ &$f_{C}$ & $2^{4}$ & $2^{2}$ &  $21$\\
			7  & $x^{22q+1}+x^{18q+5}+x^{8q+15}+x^{4q+19}+x^{23}$& $\gcd(m,11)=1$& -- & -- & --  & --  \\
			8*   & $x^{24q+1}+x^{17q+8}+x^{9q+16}+x^{8q+17}+x^{25}$& $m$ is odd & $f_{B}$  &$2^{3}$ & $2^{4}$& $25$ \\
			9*   &$x^{34q+1}+x^{33q+2}+x^{3q+32}+x^{q+34}+x^{35}$&$m$ is odd and $m\not\equiv 3\pmod{6}$&$f_{A}$& $2^{5}$ & $2$ & $35$\\
			10*   &$x^{36q+1}+x^{33q+4}+x^{32q+5}+x^{5q+32}+x^{37}$&$m\not\equiv 0\pmod{18}$& $f_{B}$&$2^{2}$ &$2^{5}$  & 37 \\
			11*    &$x^{40q+1}+x^{33q+8}+x^{9q+32}+x^{q+40}+x^{41}$& $m\not\equiv 0\pmod{10}$&  $f_{A}$  &  $2^{5}$       & $2^{3}$ & 41  \\
			12*   & $x^{48q+1}+x^{33q+16}+x^{32q+17}+x^{17q+32}+x^{49}$&$\gcd(m,3)=1$& $f_{B}$ & $2^{5}$ & $2^{4}$ &49 \\
			13*   & $x^{66q+1}+x^{65q+2}+x^{3q+64}+x^{2q+65}+x^{67}$&$m\not\equiv 0\pmod{66}$ & $f_{B}$  & $2$ &$2^{6}$  & 67 \\
			14*  & $x^{69q}+x^{64q+5}+x^{4q+65}+x^{q+68}+x^{69}$&$m\not\equiv 11\pmod{22}$& $f_{C}$& $2^{6}$ &$2^{2}$ & 69  \\
			15*  & $x^{72q+1}+x^{65q+8}+x^{9q+64}+x^{8q+65}+x^{73}$&$m \not\equiv 0\pmod{9}$& $f_{B}$& $2^{3}$ &$2^{6}$ & 73  \\
			16*  & $x^{81q}+x^{64q+17}+x^{16q+65}+x^{q+80}+x^{81}$&$m$ is odd & $f_{C}$& $2^{6}$ &$2^{4}$ & 81  \\
			17*  & $x^{96q+1}+x^{65q+32}+x^{33q+64}+x^{32q+65}+x^{97}$&$m\not\equiv 0\pmod{24}$& $f_{B}$& $2^{6}$ &$2^{5}$ & 97  \\
			\hline
	\end{tabular}}
	\caption{Permutation pentanomials over $\mathbb{F}_{2^{2m}}$ satisfying (\ref{1:fx}) and (\ref{1:gcd}) with $4 \le t<100$}
	\label{1:tab}
\end{table}
While the techniques employed in \cite{zhang2024more} to study permutation properties of $f(x)$ of the form (\ref{1:fx}) were well-understood and have been successfully applied for many similar problems in the literature, to find the exact conditions under which these 17 pentanomials are permutations over $\f$ is not a simple matter. The main technical difficulty is to study when the corresponding fractional polynomials $g(x)$ of certain degree permute $\mu_{q+1}$, the set of $(q+1)$-th roots of unity in $\f$. This is usually quite non-trivial and requires a long and hard computation for each individual family, and generally speaking, the larger the degree of $g(x)$ is, the more complex the computation becomes. This is indeed the case in \cite{zhang2024more}, and the authors managed to treat each of these 17 pentanomials individually and successfully in order to conclude that they are permutation under certain conditions. Interested readers may review the proofs of say Theorem 3.1, Lemma 3.2, Theorem 3.3, Theorem 3.4 etc. of \cite{zhang2024more} to get a taste of the techniques involved in the paper.

We first observe that among these 17 pentanomials, actually 14 of them have very simple explanation as to why they are permutations over $\f$ under certain conditions. We have listed these 14 families with a $*$-sign after their order in Table \ref{1:tab}. For example, in what below we give a new and simple explanation as to why the most complicated case, family 17 is a permutation over $\f$ when $m \not \equiv 0 \pmod{24}$:
\begin{Thm} \label{1:f17}
	Let $q=2^m$. Then $f(x)=x^{96q+1}+x^{65q+32}+x^{33q+64}+x^{32q+65}+x^{97}$ permutes $\f$ if and only if $m \not \equiv 0 \pmod{24}$.
\end{Thm}
\begin{proof} Let $\omega,\omega^2 \in \mathbb{F}_4$ be the two distinct roots of $Q(x)=x^2+x+1$.
	
	\noindent {\bf Case 1: $m$ is even.} Noting that $97=2^6+2^5+1$, we can easily check that
	\begin{eqnarray} \label{1:f97} f(x)=L_1 \circ x^{97} \circ L_2,\end{eqnarray}
	where $L_1,L_2: \f \to \f$ are given by
	\begin{eqnarray} \label{1:f17-L12} L_1(x)=\omega^2 x^q+x, \quad L_2(x)=x^q+\omega x.\end{eqnarray}
	It is also easy to check that $L_1$ and $L_2$ permute $\f$ since $m$ is even. For example, $L_1$ permutes $\f$ because $\det\left(\begin{array}{cc}
		\omega^2 & 1\\
		1 & \omega^{2q}
	\end{array}\right)=\omega^2\omega^{2q}+1=\omega^4+1=\omega^2 \ne 0.$
	
	In view of (\ref{1:f97}), we see that $f(x)$ permutes $\f$ if and only if $x^{97}$ permutes $\f$, which is equivalent to
	\begin{eqnarray} \label{1:97} \gcd\left(97,q^2-1\right)=1,
	\end{eqnarray}
	which is further equivalent to the condition that $m \not \equiv 0 \pmod{24}$, because $97$ is a prime, the order of $2$ modulo $97$ is $48$ ($2^{48} \equiv 1 \pmod{97}$ and $48$ is the least positive integer satisfying this property).
	
	\noindent {\bf Case 2: $m$ is odd.} In this case, we can easily check that
	\begin{eqnarray}
		\label{1:f17-2} f(x)=L_{1} \circ (x^{97},y^{97}) \circ L_{2},
	\end{eqnarray}
	where $L_{1}:\h^{2}\to \f$ and $L_{2}:\f \to \h^{2}$ are linear maps defined  by \begin{eqnarray} \label{1:f17-L12-2} L_{1}(x,y)=\omega x+\omega^{2} y, \quad L_{2}(x)=\left(\omega^{2}x^{q}+\omega x, \, \omega x^{q}+\omega^{2} x\right).\end{eqnarray}
	Clearly $L_{1}$ and $L_{2}$ define bijections, since $m$ is odd, $\omega\in\u$ and $\det\left(\begin{array}{cc}
		\omega^{2} & \omega  \\
		\omega & \omega^{2}
	\end{array}\right)=\omega+\omega^2=1.$
	
	In view of (\ref{1:f17-2}), we see that $f(x)$ permutes $\f$ if and only if $\left(x^{97},y^{97}\right)$ permutes $\h^{2}$, which is equivalent to the condition \begin{eqnarray}
		\label{2:f97(a)} \gcd(97, q-1)=1,
	\end{eqnarray}
	which indeed holds since $q=2^m \not\equiv 1\pmod{97}$ as $m$ is odd in this case.
	
\end{proof}

It turns out 14 of these 17 pentanomials listed in Table \ref{1:tab} with $*$-signs are all of this nature. Moreover, based on this observation, we can obtain three general classes $f_A(x),f_B(x),f_C(x)$ of permutation pentanomials satisfying (\ref{1:fx}) that these 14 pentanomials are special cases. We state these general results below.

\begin{Thm} \label{1:thmA}
	Let $m,i,j$ be positive integers and let $q=2^m, Q_{1}=2^{i}, Q_{2}=2^{j}$. The polynomial $f_A(x) \in \f[x]$ given by
	\begin{equation} \label{1:fa}
		f_{A}(x)=x^{q(Q_{1}+Q_{2})+1}+x^{q(Q_{1}+1)+Q_{2}}+x^{q(Q_{2}+1)+Q_{1}}+
		x^{Q_{1}+Q_{2}+q}+x^{Q_{1}+Q_{2}+1}
	\end{equation}
	is a permutation over $\f$ if and only if one of the following holds:
	\begin{enumerate}
		\item [i).] $m$ is odd, $ij\equiv 1\pmod2$, and $\gcd\left(Q_{1}+Q_{2}+1,q-1\right)=1$. Moreover, in this case, $f_{A}(x)$ is linear equivalent to the function  $\left(u^{Q_{1}+Q_{2}+1},v^{Q_{1}+Q_{2}+1}\right):\h^{2}\rightarrow\h^{2}$;
		
		\item[ii).] $m$ is even, $\gcd\left(Q_{1}+Q_{2}+1,q-1\right)=\gcd\left(Q_{1}+Q_{2}+1-2 r_A,q+1\right)=1$, where the integer $r_A$ is given by
		\begin{eqnarray} \label{1:ma}
			r_A=\left\{\begin{array}{lll}
				0&:& i \equiv 1, j \equiv 0 \pmod{2};\\
				1&:& i \equiv j \equiv 0 \pmod{2};\\
				Q_1&:& i \equiv 0, j \equiv 1 \pmod{2}, Q_1\le Q_2;\\
				Q_2+1&:& i \equiv 0,j \equiv 1 \pmod{2}, Q_1 > Q_2;\\
				Q_1+1&:& i \equiv 1,j \equiv 0 \pmod{2}, Q_1 < Q_2;\\
				Q_2&:& i \equiv 1,j \equiv 0 \pmod{2}, Q_1 \ge Q_2.
			\end{array}\right.
		\end{eqnarray}
		Moreover, in this case, $f_A(x)$ is linear equivalent to the function $x^{Q_1+Q_2+1+r_A(q-1)}: \f \rightarrow \f$.
	\end{enumerate}
\end{Thm}

\begin{Thm} \label{1:thmB}
	Let $m,i,j$ be positive integers and let $q=2^m, Q_{1}=2^{i}, Q_{2}=2^{j}$. The polynomial $f_{B}(x)\in \f[x]$ given by
	\begin{equation} \label{1:fb}
		f_{B}(x)=x^{q(Q_{1}+Q_{2})+1}+x^{q(Q_{1}+1)+Q_{2}}+x^{qQ_{1}+Q_{2}+1}+
		x^{q(Q_{2}+1)+Q_{1}}+x^{Q_{1}+Q_{2}+1}
	\end{equation}
	is a permutation over $\f$ if and only if one of the following holds:
	\begin{enumerate}
		\item [i).] $m$ is odd, $i(j+1)\equiv 1\pmod2$, and $\gcd\left(Q_{1}+Q_{2}+1,q-1\right)=1$. Moreover, in this case, $f_{B}(x)$ is linear equivalent to  the function $\left(u^{Q_{1}+Q_{2}+1},v^{Q_{1}+Q_{2}+1}\right):\h^{2}\rightarrow\h^{2}$;
		
		\item[ii).] $m$ is even, $\gcd\left(Q_{1}+Q_{2}+1,q-1\right)=\gcd\left(Q_{1}+Q_{2}+1-2r_B,q+1\right)=1$ where the integer $r_B$ is given by
		\begin{eqnarray}
			r_B=\left\{\begin{array}{lll}\label{1:mb}
				0&:& i \equiv 1, j \equiv 0 \pmod{2};\\
				1&:& i \equiv 0, j \equiv 1 \pmod{2};\\
				Q_1&:& i \equiv j \equiv 0 \pmod{2}, Q_1 \le Q_2;\\
				Q_2+1&:& i \equiv j \equiv 0 \pmod{2}, Q_1 > Q_2;\\
				Q_1+1&:& i \equiv j \equiv 1 \pmod{2}, Q_1 < Q_2;\\
				Q_2&:& i \equiv j \equiv 1 \pmod{2}, Q_1 \ge Q_2.
			\end{array}\right.
		\end{eqnarray}
		Moreover, in this case, $f_B(x)$ is linear equivalent to the function $x^{Q_1+Q_2+1+r_B(q-1)}: \f \rightarrow \f$.
	\end{enumerate}
\end{Thm}

\begin{Thm} \label{1:thmC}
	Let $m,i,j$ be positive integers and let $q=2^m, Q_{1}=2^{i}, Q_{2}=2^{j}$. The polynomial $f_{C}(x)\in \f[x]$ given by
	\begin{equation} \label{1:fc}
		f_{C}(x)= x^{q(Q_{1}+Q_{2}+1)}+x^{qQ_{1}+Q_{2}+1}+x^{qQ_{2}+Q_{1}+1}+x^{q+Q_{1}+Q_{2}}+x^{Q_{1}+Q_{2}+1} 		 \end{equation}
	is a permutation over $\f$ if and only if one of the following holds:
	\begin{enumerate}
		\item[i).] $m$ is odd, $(1+i)(j+1)\equiv 1\pmod 2$, and $\gcd(Q_{1}+Q_{2}+1,q-1)=1$. Moreover, in this case, $f_{C}(x)$ is linear equivalent to  the function $(u^{Q_{1}+Q_{2}+1},v^{Q_{1}+Q_{2}+1}):\h^{2}\rightarrow\h^{2}$.
		
		\item[ii).] $m$ is even, $\gcd\left(Q_{1}+Q_{2}+1,q-1\right)=\gcd\left(Q_{1}+Q_{2}+1-2r_C,q+1\right)=1$ where the integer $r_C$ is given by
		\begin{eqnarray} \label{1:mc}
			r_C=\left\{\begin{array}{lll}
				1&:& i \equiv j \equiv 1 \pmod{2};\\
				Q_1&:& i \equiv 1, j \equiv 0 \pmod{2}, Q_1 \le Q_2;\\
				Q_2+1&:& i \equiv 1,j \equiv 0 \pmod{2},Q_1 > Q_2;\\
				Q_1+1&:& i \equiv 0, j \equiv 1 \pmod{2}, Q_1 < Q_2;\\
				Q_2&:& i \equiv 0, j \equiv 1 \pmod{2}, Q_1 \ge Q_2.
			\end{array}\right.
		\end{eqnarray}
		Moreover, in this case, $f_C(x)$ is linear equivalent to the function $x^{Q_1+Q_2+1+r_C(q-1)}: \f \rightarrow \f$.
	\end{enumerate}
\end{Thm}

In Table \ref{1:tab}, we list explicitly the individual classes $f_A,f_B$ or $f_C$ that these 14 families with $*$-signs belong to, we also list the values $Q_1,Q_2$ and $Q_1+Q_2+1$.

It can be checked that \cite[Theorem 4.2]{xu2018some} is a particular instance of Theorem \ref{1:thmB} (ii) (with $Q_1=4,Q_2=2,t=7$ and $r_B=1$), and \cite[Example 4]{deng2019more} is a particular instance of Theorem \ref{1:thmB} (ii) (with $Q_1=2^k,Q_2=2,t=2^k+3$ and $r_B=0$ or $3$).

It is possible to prove Theorems \ref{1:thmA}--\ref{1:thmC} directly in a way similar to that of Theorem \ref{1:f17}, though it is conceivable that the explicit expressions of $L_1,L_2$ (see (\ref{1:f17-L12}) and (\ref{1:f17-L12-2})) may be a little complicated given the complex nature of Theorems \ref{1:thmA}--\ref{1:thmC}. We will adopt a more conventional strategy: for $\bullet \in \{A,B,C\}$, the secret of $f_\bullet(x)$ lies in the two properties Eq (\ref{3:gxeq1}) and Eq (\ref{3:gxeq2}), from which the corresponding rational function $g_\bullet(x)$ can be shown to has a very simple form (see Lemma \ref{3:clag}), from which the linear equivalence of $f_\bullet(x)$ arises naturally (see Lemma \ref{3:clag2}), and this will complete the proofs of Theorems \ref{1:thmA}--\ref{1:thmC}.

We remark that the idea of treating the rational function $g_\bullet(x)$ on $\mu_{q+1}$ in such a way was originated from the paper \cite{ding2023determination}. The idea of linear equivalence was also hinted in \cite{ding2023determination}. We acknowledge that the proofs of Theorems \ref{1:thmA}--\ref{1:thmC} are very simple applications of ideas from \cite{ding2023determination}.

This paper is organized as follows: in Section \ref{2:pre} we recall some basic terminologies that will be used throughout the paper; in Section \ref{3:chap3} we study basic properties of the numerator and the denominator of the rational function $g_\bullet(x)$; in Section \ref{4:chap4} we use geometric properties of $g_\bullet(x)$ to prove Lemmas \ref{3:clag} and \ref{3:clag2}, from which we prove Theorems \ref{1:thmA}--\ref{1:thmC} directly. In Section \ref{fin:cond} we conclude the paper.

\section{Preliminaries} \label{2:pre}

In this section, we shall recall some basic terminologies that will be used to prove the main results. We begin with some notations:	
\begin{itemize}[itemsep=0em]
	\item For $x\in\f$, denote $\overline{x}=x^{q};$
	\item $\y$ is the algebraic closure of $\h$;
	\item the unit circle of $\f$ is denoted as the set $\u=\left\{x\in\f \mid x^{q+1}=\overline{x}x=1\right\};$
	\item $\z=\y\cup\{\infty\};$
\end{itemize}

\subsection{Geometric properties of rational functions}

Here we introduce the concepts of branch points and ramification, which are the main tools to describe the geometric properties of the accompanying rational function $g(x)$ later on. While all the details were already given in \cite{ding2023determination}, for simplicity we follow the presentation from \cite{Terry2024}.

For a non-constant rational function $G(x) \in \overline{\F}_q(x)$, write $G(x)=N(x)/D(x)$ where $N(x), D(x) \in \overline{\F}_q[x]$ with $\gcd(N(x),D(x))=1$. For any $\alpha \in \overline{\F}_q$, define $H_\alpha(x) \in \overline{\F}_q[x]$ as
$$H_\alpha(x):=\begin{cases}
	N(x)-G(\alpha)D(x) &(G(\alpha) \in \overline{\F}_q),\\
	D(x) &(G(\alpha)=\infty).
\end{cases}$$
It is clear that $H_\alpha(\alpha)=0$ for all $\alpha \in \overline{\F}_q$. The \emph{ramification index} $e_G(\alpha)$ of $\alpha$ is then its multiplicity as a root of $H_\alpha(x)$. The ramification index of $\infty$ is defined as $e_G(\infty):=e_{G_1}(0)$, where $G_1(x):=G(\frac{1}{x})$. Given any $\beta \in \P^1(\overline{\F}_q)$, the \emph{ramification multiset} $E_G(\beta)$ of $G(x)$ over $\beta$ is the multiset of ramification indices $e_G(\alpha)$ for any $\alpha \in G^{-1}(\beta)$, that is, $E_G(\beta):=\left[e_G(\alpha): \alpha \in G^{-1}(\beta)\right]$. In particular, the elements in $E_G(\beta)$ are positive integers whose sum is $\deg G$. We call $\alpha \in \P^1(\overline{\F}_q)$ a \emph{ramification point} (or \emph{critical point}) of $G(x)$ if $e_G(\alpha) > 1$, and its corresponding image $G(\alpha)$ a \emph{branch point} (or \emph{critical value}) of $G$. Hence a point $\beta \in \P^1(\overline{\F}_q)$ is a branch point of $G$ if and only if $E_G(\beta) \neq [1^{\deg G}]$, where $[m^n]$ denotes the multiset consisting of $n$ copies of $m$, or equivalently, $\sharp G^{-1}(\beta) < \deg G$.

In the next result, for any polynomial $P(x)\in\y[x]$ we write $P'(x)$ for the (formal) derivative of $P(x)$.

\begin{Lem}\cite[Lemma 2.12]{Ding2024}
	\label{branchpoint, dervative}
	Let $h(x)=\frac{N(x)}{D(x)} \in\y(x)$ be a non-constant rational function where $N,D \in \y[x]$ with $\gcd\left(N(x),D(x)\right)=1$. If $\alpha \in \y$ is a ramification point of $h(x)$, then
	$$N'(\alpha)D(\a)-N(\a)D'(\a)=0.$$
\end{Lem} 	


\begin{Lem}\cite[Lemma 2.1]{zieve2013permutation}
	\label{u2u}
	Let $q$ be a prime power, and let $\rho(x) \in\y[x]$ be a degree one rational function. Then $\rho$ induces a bijection on $\mu_{q+1}$ if and only if $\rho(x)$ equals either
	\begin{enumerate}
		\item $\b x$ or $\b/x$ with $\b \in\mu_{q+1},$ or
		\item $(x-\gamma^{q}\beta)/(\gamma x-\b)$ with $\b\in\mu_{q+1}$ and $\gamma \in\mathbb{F}_{q^{2}}\backslash \mu_{q+1}$.
	\end{enumerate}
\end{Lem}

\begin{Lem}\cite[Lemma 3.1]{zieve2013permutation}
	\label{u2P}
	Let $q$ be a prime power and $l\in\y(x)$ be a degree-one rational function. Then $l(x)$ induces a bijection from $\mu_{q+1}$ to $\pp$ if and only if $l(x)=(\delta x- \beta \delta^{q})/(x-\beta)$ with $\beta\in \mu_{q+1}$ and $\delta\in\mathbb{F}_{q^{2}}\backslash\mathbb{F}_{q}.$
\end{Lem}

\subsection{Linear equivalence}

Let us quickly review linear equivalence, EA-equivalence and CCZ-equivalence of boolean functions, the three important equivalence relations of functions over the finite field $\mathbb{F}_{p^n}$. For the presentation, we focus only on the linear equivalence. Interested readers may refer to \cite{shi2024ccz, budaghyan2006new, carlet1998codes} for more information about the other two equivalence relations.

\begin{definition} A function $F: \F_{p^n} \to \F_{p^n}$ is called
	\begin{itemize}
		\item \emph{linear} if $F(\alpha+\beta)=F(\alpha)+F(\beta)$ for any $\alpha, \beta \in \F_{p^n}$;
		
		\item \emph{affine} if $F$ is a sum of a linear function and a constant;
		
		\item \emph{affine permutation} (or \emph{linear permutation}) if $F$ is both affine (resp. linear) and a permutation on $\F_{p^n}$.
	\end{itemize}
	In particular, $F$ is affine if and only if $F(x)=b+\sum_{j=0}^{n-1} a_jx^{p^j}$ where $a_j,b \in \F_{p^n}$ for any $j$.
\end{definition}
\begin{definition} Two functions $F$ and $F'$ from $\F_{p^n}$ to itself are called:
	\begin{itemize}
		\item \emph{affine equivalent} (or \emph{linear equivalent}) if $F'=A_1 \circ F \circ A_2$, where the mappings $A_1,A_2$ are affine (resp. linear) permutations of $\F_{p^n}$;
		\item \emph{extended affine equivalent} (EA-equivalent) if $F'=A_1 \circ F \circ A_2+A$, where the mappings $A,A_1,A_2$ are affine, and where $A_1,A_2$ are permutations of $\F_{p^n}$;
		\item \emph{Carlet-Charpin-Zinoviev equivalent} (CCZ-equivalent) if for some affine permutation $\mathcal{L}$ of $\F_{p^n}^2$ the image of the graph of $F$ is the graph of $F'$, that is, $\mathcal{L}(G_F)=G_{F'}$, where
		\[G_F=\left\{(\a,F(\a)): \a \in \F_{p^n}\right\}, \quad G_{F'}=\left\{(\a,F'(\a)): \a \in \F_{p^n}\right\}. \]
	\end{itemize}
\end{definition}

It is obvious that linear equivalence is a particular case of affine equivalence, and affine equivalence is a particular case of EA-equivalence. If $f(x)$ and $g(x)$ are linear equivalent, then $f$ and $g$ share many properties. In particular, if $f(x)$ is a permutaiton on $\mathbb{F}_q$, then so is $g(x)$.
	
	\section{Properties of $H_\bullet(x)$ and $N_\bullet(x)$} \label{3:chap3}
Let $q=2^m, Q_{1}=2^{i},Q_{2}=2^{j}$ where $m,i,j$ are some positive integers. For any $\bullet \in \{A,B,C\}$, the polynomial $f_{\bullet}(x)$ given in (\ref{1:fa})-(\ref{1:fc}) can be written respectively as
\begin{eqnarray} \label{3:fx} f_{\bullet}(x)=x^{Q_1+Q_2+1}H_{\bullet}(x),\end{eqnarray} where
\begin{eqnarray}
	H_A(x)&=&x^{Q_{1}+Q_{2}}+x^{Q_{1}+1}+x^{Q_{2}+1}+x+1; \label{3:ha}\\
	H_B(x)&=&x^{Q_{1}+Q_{2}}+x^{Q_{1}+1}+x^{Q_{1}}+x^{Q_{2}+1}+1; \label{3:hb}\\
	H_C(x)&=&x^{Q_{1}+Q_{2}+1}+x^{Q_{1}}+x^{Q_{2}}+x+1. \label{3:hc}
\end{eqnarray}
It is well-known that $f_{\bullet}(x)$ permutes $\f$ if and only if $\gcd(Q_1+Q_2+1,q-1)=1$ and $x^{Q_1+Q_2+1}H_{\bullet}(x)^{q-1}$ permutes $\mu_{q+1}$, the set of $(q+1)$-th roots of unity in $\f$ (see for example \cite{zieve2009some}).
By a lemma introduced in \cite{ding2023constructing} (see also \cite[Lemma 2.10]{ding2023determination}), this in turn is equivalent to
\begin{itemize}
	\item[(1)] $\gcd(Q_1+Q_2+1,q-1)=1$,
	\item[(2)] $H_{\bullet}(x)$ has no roots in $\mu_{q+1}$,
	\item[(3)] $g_{\bullet}(x)$ permutes $\mu_{q+1}$,
\end{itemize}
where
\begin{eqnarray}
	 g_{A}(x)&=&\frac{N_{A}(x)}{H_{A}(x)}=\frac{x+x^{Q_{2}}+x^{Q_{1}}+x^{Q_{1}+Q_{2}}+x^{Q_{1}+Q_{2}+1}}{x^{Q_{1}+Q_{2}}+x^{Q_{1}+1}+x^{Q_{2}+1}+x+1},
	\label{3:ga} \\
	 g_{B}(x)&=&\frac{N_{B}(x)}{H_{B}(x)}=\frac{x+x^{Q_{2}}+x^{Q_{1}}+x^{Q_{2}+1}+x^{Q_{1}+Q_{2}+1}}{x^{Q_{1}+Q_{2}}+x^{Q_{1}+1}+x^{Q_{1}}+x^{Q_{2}+1}+1},
	\label{3:gb}\\
	 g_{C}(x)&=&\frac{N_{C}(x)}{H_{C}(x)}=\frac{1+x^{Q_{1}+1}+x^{Q_{2}+1}+x^{Q_{1}+Q_{2}}+x^{Q_{1}+Q_{2}+1}}{x^{Q_{1}+Q_{2}+1}+x^{Q_{1}}+x^{Q_{2}}+x+1}. \label{3:gc}
\end{eqnarray}
Here it is easy to check that that for any $\bullet \in \{A,B,C\}$,
\begin{eqnarray} \label{3:nx} N_{\bullet}(x)=x^{Q_1+Q_2+1}H_{\bullet}\left(x^{-1}\right),
\end{eqnarray}
and $g_{\bullet}(x)=\frac{N_\bullet(x)}{H_\bullet(x)}$ satisfies the following two properties
\begin{eqnarray}
	N_{\bullet}'(x)H_{\bullet}(x)+N_{\bullet}(x)H_{\bullet}'(x)&=&Q(x)^{Q_1+Q_2}, \label{3:gxeq1}\\
	Q\left(g_{\bullet}(x)\right) H_{\bullet}(x)^2&=&Q(x)^{Q_1+Q_2+1}. \label{3:gxeq2}
\end{eqnarray}
Here we denote
\begin{eqnarray} \label{3:qx} Q(x)=x^2+x+1 \in \mathbb{F}_2[x],
\end{eqnarray}
and $f'(x)$ is the formal derivative of any rational function $f(x)$. Denote by $\omega, \omega^2\in \mathbb{F}_4$ the two distinct roots of $Q(x)$. It is easy to see that $\omega$ satisfies
\begin{eqnarray} \label{3:omega} \omega \in \mathbb{F}_4, \quad \omega^3=1, \quad \omega +\omega^2=1. \end{eqnarray}
We first study properties of $H_\bullet(x)$ and $N_\bullet(x)$ for any $\bullet \in\{A,B,C\}$.


\begin{Lem}\label{3:gcd}
	For $\bullet \in \{A,B,C\}$, let $g_{\bullet}(x)=\frac{N_{\bullet}(x)}{H_{\bullet}(x)}$. Denote $C_{\bullet}(x)=\gcd(N_{\bullet}(x),H_{\bullet}(x))$ the (monic) greatest common divisor of $N_{\bullet}(x)$ and $H_{\bullet}(x)$. Then we have
	\begin{eqnarray}
		C_{\bullet}(x)&=&Q(x)^{r_{\bullet}},
		\label{3:mm} \end{eqnarray}
	where the positive integers $r_A,r_B$ and $r_C$ above are given in (\ref{1:ma}) (\ref{1:mb}) and (\ref{1:mc}) respectively.
\end{Lem}
\begin{proof}
	For any $\bullet \in \{A,B,C\}$, $H_{\bullet}(x),N_{\bullet}(x) \in \mathbb{F}_2[x]$, so $C_{\bullet}(x) \in \mathbb{F}_2[x]$; since $C_\bullet(x) \mid H_\bullet(x)$ and $C_\bullet(x) \mid N_\bullet(x)$, by Eq (\ref{3:gxeq1}), $C_{\bullet}(x) \mid Q(x)^{Q_1+Q_2}$. Finally, since $Q(x)=x^2+x+1$ is irreducible over $\mathbb{F}_2$, we have
	\begin{eqnarray} \label{3:cm} C_{\bullet}(x)=Q(x)^{r_\bullet}
	\end{eqnarray}
	for some integer $r_{\bullet} \ge 0$. To prove Lemma \ref{3:gcd}, we just need to compute these numbers $r_{\bullet}$ for $\bullet \in \{A,B,C\}$.
	
	Let $r_{\bullet}'$ be the largest nonnegative integer such that $\left(x+\omega\right)^{r_{\bullet}'} \mid H_{\bullet}(x)$. We write it as
	\begin{eqnarray} \label{3:xom} \left(x+\omega\right)^{r_{\bullet}'} \| H_{\bullet}(x).\end{eqnarray}
	Here $\omega, \omega^2$ are the two distinct roots of $Q(x)$ satisfying properties (\ref{3:omega}). Taking squares on every coefficient of polynomials on both sides of Eq (\ref{3:xom}) (this is a Frobenius action) we have $\left(x+\omega^2\right)^{r_{\bullet}'} \mid H_{\bullet}(x)$, which implies that
	\[Q(x)^{r_{\bullet}'} \mid H_{\bullet}(x).\]
	Noting that $x^2Q(x^{-1})=Q(x)$, we have
	\[Q(x)^{r_\bullet'}=x^{2r_{\bullet}'}Q\left(x^{-1}\right)^{r_\bullet'} \mid x^{Q_1+Q_2+1}H_\bullet\left(x^{-1}\right)=N_\bullet(x),\]
	that is,
	$$Q(x)^{r_\bullet'} \mid \gcd(H_\bullet(x),N_\bullet(x))=C_\bullet(x)=Q(x)^{r_\bullet}.$$
	It is also easy to see that actually $r_\bullet'$ is the largest nonnegative integer satisfying the above property, so $r_\bullet'=r_\bullet$. Thus to find $r_\bullet$, we just need to find the integer $r_\bullet'$ satisfying Eq (\ref{3:xom}).
	
	Let $i_0,j_0 \in \{0,1\}$ be such that $i \equiv i_0 \pmod{2}, j \equiv j_0 \pmod{2}$. Since $\omega$ is a root of $Q(x)$ and satisfies (\ref{3:omega}), it is easy to see that
	\[\omega^{Q_1}=\omega^{i_0+1}, \quad \omega^{Q_2}=\omega^{j_0+1}.\]
	
	\noindent {\bf Case A:} We can compute
	\begin{eqnarray*} H_A(\omega)&=&\omega^{Q_1+Q_2}+\omega^{Q_1+1}+\omega^{Q_2+1}+\omega+1\\
		&=&\omega^{i_0+j_0+2}+\omega^{i_0+2}+\omega^{j_0+2}+\omega+1\\
		&=&\omega^2 \left(\omega^{i_0}+1\right)\left(\omega^{j_0}+1\right).\end{eqnarray*}
	
	If $ij \equiv 1 \pmod{2}$, then $i_0=j_0=1$, it is easy to see that $H_A(\omega)=1$, that is, $\left(x+\omega\right) \nmid H_A(x)$, so we have $r_A=0$.
	
	Now let us assume that $ij \equiv 0 \pmod{2}$, that is, $i_0=0$ or $j_0=0$. It is easy to check that $H_A(\omega)=0$. By Eq (\ref{3:xom}), we just need to find $r_A$ such that
	\[x^{r_A} \| H_A\left(x+\omega\right). \]
	We can compute that
	\begin{eqnarray*} H_A\left(x+\omega\right)&=&\left(x+\omega\right)^{Q_1+Q_2}+\left(x+\omega\right)^{Q_1+1}+
		\left(x+\omega\right)^{Q_2+1}+\left(x+\omega\right)+1\\
		 &=&x^{Q_1+Q_2}+x^{Q_1}\left(\omega^{Q_2}+\omega\right)+x^{Q_2}\left(\omega^{Q_1}+\omega\right)
		+x^{Q_1+1}+x^{Q_2+1}+x\left(\omega^{Q_1}+\omega^{Q_2}+1\right)\\
		 &=&x^{Q_1+Q_2}+x^{Q_1}\left(\omega^{j_0+1}+\omega\right)+x^{Q_2}\left(\omega^{i_0+1}+\omega\right)
		+x^{Q_1+1}+x^{Q_2+1}+x\left(\omega^{i_0+1}+\omega^{j_0+1}+1\right). \end{eqnarray*}
	It is easy to see by checking non-zero coefficients of $H_A\left(x+\omega\right)$ that if $x^ r\|H_A(x+\omega)$, then the positive integer $r$ is exactly $r_A$ given by (\ref{1:ma}). For example, if $i \equiv j \equiv 0 \pmod{2}$, then $i_0=j_0=0$, the lowest order term is $x$, whose coefficient is  $\omega^{i_0+1}+\omega^{j_0+1}+1=1 \ne 0$ so $x \|H_A(x+\omega)$ and hence $r_A=1$ in this case; if $i \equiv 0,j \equiv 1 \pmod{2}$, then $i_0=0,j_0=1$, we have $\omega^{i_0+1}+\omega^{j_0+1}+1=0$, if $Q_1 \le Q_2$, the lowest order term of $H_A\left(x+\omega\right)$ is $x^{Q_1}$, whose coefficient is $\omega^{j_0+1}+\omega=1 \ne 0$, so $x^{Q_1} \|H_A(x+\omega)$ and hence $r_A=Q_1$; if in this case $Q_1>Q_2$ instead, then since $\omega_{i_0+1}+\omega=0$, the lowest order term is $x^{Q_2+1}$, whose coefficient is $1 \ne 0$, so $r_A=Q_2+1$ in this case. The other cases are similar. This proves (\ref{3:cm}) for {\bf Case A}.
	
	\noindent {\bf Case B:} We can obtain
	\begin{eqnarray*}
		H_{B}(\omega) &=& \omega^{i_{0}+1}(\omega^{j_{0}+1}+\omega+1)+(\omega^{j_{0}+2}+1)\\
		&=& \left(\omega^{i_{0}}+1\right)\left(\omega^{j_{0}+2}+1\right).
	\end{eqnarray*}
	If $i(j+1)\equiv 1\pmod{2}$, then $i_{0}=1,j_{0}=0$, it is easy to check that $H_{B}(\omega)=1$, that is, $(x+\omega)\nmid H_{B}(x)$, so we have $r_{B}=0$.
	
	Let us assume that $i(j+1)\equiv 0\pmod{2}$, then $i_{0}=0$ or $j_{0}=1$. It is easy to check that $H_{B}\left(\omega\right)=0$. In order to calculate $r_{B}$, we must first compute $H_{B}(x+\omega)$, where
	\begin{eqnarray*}
		 H_{B}(x+\omega)&=&(x+\omega)^{Q_{1}+Q_{2}}+(x+\omega)^{Q_{1}+1}+(x+\omega)^{Q_{2}+1}+(x+\omega)^{Q_{1}}+1\\
		&=& x^{Q_{1}+Q_{2}}+x^{Q_{1}+1}+x^{Q_{2}+1}+(\omega^{Q_{2}}+\omega^{2})x^{Q_{1}}+(\omega^{Q_{1}}+\omega)x^{Q_{2}}+(\omega^{Q_{1}}+\omega^{Q_{2}})x\\
		&=&x^{Q_1+Q_2}+x^{Q_1+1}+x^{Q_2+1}+x^{Q_1}\left(\omega^{j_0+1}+\omega^2\right)
		 +x^{Q_2}\left(\omega^{i_0+1}+\omega\right)+x\left(\omega^{i_0+1}+\omega^{j_0+1}\right).
	\end{eqnarray*}
	Using the same procedure that we used to compute $r_{A}$ in {\bf Case A}, we see that if $x^r\| H_{B}(x+\omega)$, then $r=r_B$ which was given in Eq (\ref{1:mb}). We omit the details. This proves (\ref{3:cm}) for {\bf Case B}.
	
	\noindent {\bf Case C:} In this case, we have
	\begin{eqnarray*}
		H_{C}(\omega)&=&\omega^{i_{0}}(\omega^{j_{0}}+\omega)+(\omega^{j_{0}+1}+\omega+1)\\
		&=&\left(\omega^{i_0}+\omega\right)\left(\omega^{j_0}+\omega\right).
	\end{eqnarray*}
	If $(1+i)(1+j)\equiv 1\pmod{2}$, then $i_{0}=j_{0}=0$, and clearly $H_{C}(\omega)=\omega \ne 0$, that is, $(x+\omega)\nmid H_{C}(x)$, thus we have $r_{C}=0$ in this case.
	
	If $((1+j)(1+i))\equiv 0\pmod{2}$, then $i_{0}=1$ or $j_{0}=1$, it is easy to see that $H_{C}(\omega)=0$. In order to calculate $r_{C}$, we must first compute $H_{C}(x+\omega)$, where
	\begin{equation*}
		\begin{split}
			&H_C\left(x+\omega\right)\\
			&=(x+\omega)^{Q_{1}+Q_{2}+1}+(x+\omega)^{Q_{1}}+(x+\omega)^{Q_{2}}+(x+\omega)+1\\
			&= x^{Q_{1}+Q_{2}}\left(x+\omega\right)+\omega^{Q_{2}}x^{Q_{1}}+\omega^{Q_{1}}x^{Q_{2}}+(\omega^{Q_{2}+1}+1)x^{Q_{1}}+(\omega^{Q_{1}+1}+1)x^{Q_{2}}+(\omega^{Q_{1}+Q_{2}}+1)x\\
			 &=x^{Q_{1}+Q_{2}}(x+\omega)+\omega^{j_{0}+1}x^{Q_{1}+1}+\omega^{i_{0}+1}x^{Q_{2}+1}+(\omega^{j_{0}+2}+1)x^{Q_{1}}+(\omega^{i_{0}+2}+1)x^{Q_{2}}+(\omega^{i_{0}+j_{0}+2}+1)x.
		\end{split}
	\end{equation*}
	Using the same procedure that we used to compute $r_{A}$ in {\bf Case A}, we see that if $x^r\| H_{C}(x+\omega)$, then $r=r_C$ which was given in Eq (\ref{1:mc}). We omit the details. This proves (\ref{3:cm}) for {\bf Case C}.
	
\end{proof}

\begin{Lem}\label{3:hxo}
	For $\bullet \in \{A,B,C\}$, let $g_{\bullet}(x)=\frac{N_{\bullet}(x)}{H_{\bullet}(x)}$. \begin{enumerate}[itemsep=0em]
		\item [i).] $H_A(x)$ has no roots in $\mu_{q+1}$ if and only if $ij\equiv 1\pmod{2}$ or $m$ is even;
		\item [ii).] $H_B(x)$ has no roots in $\mu_{q+1}$ if and only if $i(j+1)\equiv 1\pmod{2}$ or $m$ is even;
		\item [iii).] $H_C(x)$ has no roots in $\mu_{q+1}$ if and only if $(1+i)(1+j)\equiv 1\pmod{2}$ or $m$ is even.
	\end{enumerate}
\end{Lem}
\begin{proof}
	
	Let $\alpha \in \mu_{q+1}$ be such that $H_\bullet(\alpha)=0$. Since $\alpha^q=\alpha^{-1}$, we have $ H_\bullet(\alpha)^q=H_\bullet\left(\alpha^{-1}\right)=0$, thus  $N_\bullet(\alpha)=\alpha^{Q_1+Q_2+1}H_\bullet\left(\alpha^{-1}\right)=0$. So we have $ \left(x+\alpha\right)|H_A(x)$ and $\left(x+\alpha\right)|N_A(x)$, which implies that
	\begin{eqnarray} \label{3:alpha} \left(x+\alpha\right)|\gcd\left(H_A(x),N_A(x)\right)=C_A(x)=Q(x)^{r_\bullet}.\end{eqnarray}
	So $r_\bullet>0$. The two distinct roots of $Q(x)$ are $\omega,\omega^2 \in \mathbb{F}_4$. Eq (\ref{3:alpha}) now implies that $\alpha \in\{\omega,\omega^2\}$. On the other hand, it is easy to check that $\omega \in \mu_{q+1}$ if and only if $m$ is odd. This shows that $H_\bullet(\alpha)=0$ for some $\alpha \in \mu_{q+1}$ if and only if $m$ is odd and $r_\bullet>0$. Checking the values of $r_\bullet$ for $\bullet \in \{A,B,C\}$ from (\ref{1:ma}) (\ref{1:mb}) and (\ref{1:mc}) respectively proves Lemma \ref{3:hxo}.
\end{proof}

	\section{Proof of Theorems \ref{1:thmA}--\ref{1:thmC}} \label{4:chap4}
In this section, we first study geometric properties of $g_\bullet(x)$, from which Theorems \ref{1:thmA}--\ref{1:thmC} can be proved easily.

\begin{Lem}
	\label{3:clag}
	For $\bullet \in \{A,B,C\}$, let $g_{\bullet}(x)=\frac{N_{\bullet}(x)}{H_{\bullet}(x)}$. Assume that $\gcd\left(N_\bullet(x),H_\bullet(x)\right)=Q(x)^{r_{\bullet}}$ and $H_\bullet(x)$ has no roots in $\mu_{q+1}$. Then $g_\bullet(x)$ has the following classification.
	\begin{enumerate}[itemsep=0em]
		\item If $m$ is odd, then $r_\bullet=0$ and $g_\bullet(x)=\eta^{-1}\circ x^{Q_{1}+Q_{2}+1}\circ \sigma$ for some degree-one $\eta,\sigma\in\f(x)$ both of which map $\u$ onto $\pp$;
		\item If $m$ is even, then $g_\bullet(x)=\rho^{-1}\circ x^{Q_{1}+Q_{2}+1-2r_\bullet}\circ\sigma$ for some degree-one $\rho,\sigma\in\f(x)$ both of which permute $\u$.
	\end{enumerate}
\end{Lem}

\begin{proof}
	For $\bullet \in \{A,B,C\}$, denote
	\begin{eqnarray*}
		H_\bullet(x)=\widetilde{H}_{\bullet}(x)Q(x)^{r_\bullet} , \quad N_\bullet(x)= \widetilde{N}_{\bullet}(x)Q(x)^{r_\bullet}.
	\end{eqnarray*}
	Then
	\[\gcd\left(\widetilde{H}_\bullet(x),\widetilde{N}_\bullet(x)\right)=1,\]
	and
	 \[g_\bullet(x)=\frac{N_\bullet(x)}{H_\bullet(x)}=\frac{\widetilde{N}_\bullet(x)}{\widetilde{H}_\bullet(x)}.\]
	Eq (\ref{3:gxeq1}) can be written as
	\begin{equation}
		\label{3:ram}
		g_\bullet'(x)\widetilde{H}_\bullet(x)^{2}=Q(x)^{Q_{1}+Q_{2}-2r_\bullet},
	\end{equation}
	and Eq (\ref{3:gxeq2}) can be written as
	\begin{equation}
		\label{3:branch}
		 Q\left(g_\bullet(x)\right)\widetilde{H}_\bullet(x)^{2}=Q(x)^{Q_{1}+Q_{2}+1-2r_\bullet}.
	\end{equation}
	Eq (\ref{3:ram}) means that the set of ramification points of $g_\bullet(x)$ is $\left\{\omega,\omega^2\right\}$, where $\omega,\omega^2$ are the two distinct roots of $Q(x)$ and since \[Q\left(g_\bullet(x)\right)=\left(g_\bullet(x)+\omega\right)\left(g_\bullet(x)+\omega^2\right),\]
	Eq (\ref{3:branch}) means that the set of branch points of $g_\bullet(x)$ is $\left\{g_\bullet\left(\omega\right),g_\bullet\left(\omega^2\right)\right\}=\left\{\omega,\omega^2\right\}$ and the ramification indices of $g_\bullet(x)$ are given by
	 \[E_{g_\bullet}\left(\omega\right)=E_{g_\bullet}\left(\omega^2\right)=\left[Q_1+Q_2+1-2r_\bullet\right].\]
	Noting $\deg g_{\bullet}(x)=Q_1+Q_2+1-2r_\bullet$, we see that Lemma \ref{3:clag} is a direct consequence of \cite[Lemma 5.1]{ding2023determination} (with a very small variation). For the sake of completeness, we provide a detailed proof here, because the proof is quite simple for our particular case.
	
	For simplicity let us denote
	\[\alpha_k:=\omega^k, \quad \beta_k:=g_\bullet\left(\omega^k\right), \quad k=1,2.\]
	So the ramification points and the corresponding branch points of $g_\bullet(x)$ are $\alpha_k$ ($k=1,2$) and $\beta_k$ ($k=1,2$) respectively.
	
	\noindent \textbf{Case 1: $m$ is odd.}
	
	\noindent In this case $r_\bullet=0$ and $\omega \in \u$. So $\a_{k},\b_k \in \u$ for $k=1,2$. Picking distinct elements $\delta_{1},\delta_{2}\in\f\setminus\h$ such that
	\[ \delta_{1}^{q-1}=\a_{2}/\a_{1}, \quad \delta_{2}^{q-1}=\b_{2}/\b_{1},\]
	we can define two degree-one rational functions
	$$\sigma(x):=\frac{\delta_{1}\left(x+\a_{2}\right)}{x+\a_{1}},\quad \rho(x):=\frac{\delta_{2}\left(x+\b_{2}\right)}{x+\b_{1}}.$$
	By Lemma \ref{u2P}, both $\sigma(x)$ and $\rho(x)$ induce bijections from $\mu_{q+1}$ to $\pp$. Let us define
	\begin{eqnarray*} h(x):=\rho \circ g_\bullet(x)\circ \sigma^{-1}.\end{eqnarray*}
	Since $\sigma(\a_{1})=\rho(\b_{1})=\infty$ and $\sigma(\a_{2})=\rho(\b_{2})=0$, we have $h(0)=0$ and $h(\infty)=\infty$, and
	\[E_{h}\left(0\right)=E_{h}\left(\infty\right)=\left[Q_1+Q_2+1\right].\]
	Since $\deg h=Q_1+Q_2+1$, we have that $h(x)=c x^{Q_{1}+Q_{2}+1}$ for some $c\in\y^{\ast}$. Here actually $c\in\g$ as $h\in \mathbb{F}_q(x)$ by \cite[Lemma 2.11]{ding2023determination}. Letting $\eta(x)=c^{-1}\rho(x)$, $\eta$ is a bijection from $\mu_{q+1}$ to $\pp$, thus we have
	\[g_\bullet(x)=\eta^{-1} \circ x^{Q_{1}+Q_{2}+1} \circ \sigma,\]
	as claimed by Case 1 of Lemma \ref{3:clag}.

	\noindent \textbf{Case 2: $m$ is even.}
	\noindent In this case, $\omega \in \g$. So $\a_{k},\b_k \in \g$ for $k=1,2$. We define two degree-one rational functions
	$$\rho(x):=\frac{\b_{1}\left(x+\b_2\right)}{x+\b_{1}}, \quad \sigma(x):= \frac{\a_{1}\left(x+\a_2\right)}{x+\a_{1}}.$$
	By Lemma \ref{u2u}, both $\rho(x),\sigma(x)$ induce permutations on $\mu_{q+1}$. Let us define $$h(x)=\rho \circ g_\bullet(x) \circ \sigma^{-1}.$$
	Since $\rho(\b_{1})=\sigma(\a_{1})=\infty$ and $\rho(\b_{2})=\sigma(\a_{2})=0$, we have $h(0)=0$ and $h(\infty)=\infty$, and
	\[E_{h}\left(0\right)=E_{h}\left(\infty\right)=\left[Q_1+Q_2+1-2r_\bullet\right].\]
	Since $\deg h=Q_1+Q_2+1-2r_\bullet$, again we have $h(x)=c x^{Q_{1}+Q_{2}+1-2r_\bullet}$ for some
	$c\in\g$ by using \cite[Lemma 2.11]{ding2023determination}. Letting $\eta(x)=c^{-1}\rho(x)$, thus we have
	\[g_\bullet(x)=\eta^{-1} \circ x^{Q_{1}+Q_{2}+1-2r_\bullet} \circ \sigma,\]
	as claimed by Case 2 of Lemma \ref{3:clag}.
\end{proof}
Armed with Lemma \ref{3:clag}, we show that a linear equivalence of $f_\bullet(x)$ given in (\ref{3:fx}) can be derived easily. The idea was first hinted in the proof of \cite[Theorem 1.2]{ding2023determination}. It was used in a similar setting in \cite[Lemma 15]{Terry2024}. For the sake of completeness, we provide a detailed proof here.

\begin{Lem}
	\label{3:clag2}
	For $\bullet \in \{A,B,C\}$, let $f_\bullet(x)$ be defined in (\ref{3:fx}) and let $g_{\bullet}(x)=\frac{N_{\bullet}(x)}{H_{\bullet}(x)}$. Assume that $\gcd\left(N_\bullet(x),H_\bullet(x)\right)=Q(x)^{r_{\bullet}}$ and $H_\bullet(x)$ has no roots in $\mu_{q+1}$.
	\begin{enumerate}
		\item If $m$ is odd, then $f_\bullet(x)$ is linear equivalent to the bivariate function
		$$P_{1}(u,v)=\left(u^{Q_{1}+Q_{2}+1}, v^{Q_{1}+Q_{2}+1}\right):\h^{2}\rightarrow\h^{2};$$
		\item If $m$ is even, then $f_\bullet(x)$ is linear equivalent to the function
		$$P_{2}(x)=x^{Q_{1}+Q_{2}+1+r_\bullet(q-1)}:\f\rightarrow\f.$$
	\end{enumerate}
\end{Lem}

\begin{proof}
	\noindent {\bf Case 1: $t$ is odd.} From 1 of Lemma \ref{3:clag}, there are degree-one rational functions $\eta(x), \sigma(x)$ over $\f$ defined from $\u$ to $\pp$ such that $g_\bullet(x)=\eta^{-1}\circ x^{Q_{1}+Q_{2}+1}\circ \sigma$. By Lemma \ref{u2P}, we may write $\eta$ and $\sigma$ as
	\[\eta(x)=\frac{\a_{2}x+\b_{2}\overline{\a}_{2}}{x+\b_{2}},\quad \sigma(x)= \frac{\a_{1}x+\b_{1}\overline{\a}_{1}}{x+\b_{1}},\]
	where $\a_{1},\a_{2}\in\f\setminus\h$ and $\b_{1},\b_{2}\in\u$, and $\overline{\a}:=\a^q$ for any $\a \in \f$. Thus
	\[g_\bullet (x)=\frac{\b_{2}(x+\overline{\a}_{2})}{x+\a_{2}}\circ x^{Q_{1}+Q_{2}+1}\circ \frac{\a_{1}x+\b_{1}\overline{\a}_{1}}{x+\b_{1}},\]
	We may expand the above expression to obtain
	\begin{eqnarray} \label{3:1gx} g_\bullet(x)=\frac{\b_2\left(\left(\a_{1}x+\b_{1}\overline{\a}_{1}\right)^{Q_{1}+Q_{2}+1}+\overline{\a}_{2}\left(x+\b_{1}\right)^{Q_{1}+Q_{2}+1}\right)}{\left(\a_{1}x+\b_{1}\overline{\a}_{1}\right)^{Q_{1}+Q_{2}+1}+\a_{2}\left(x+\b_{1}\right)^{Q_{1}+Q_{2}+1}}.\end{eqnarray}
	Noting that in this case $g_\bullet(x)=\frac{N_\bullet(x)}{H_\bullet(x)}$ with $\gcd\left(N_\bullet(x),H_\bullet(x)\right)=1$. Since $\deg N_\bullet(x)=Q_1+Q_2+1$ for any $\bullet \in \{A,B,C\}$, Eq (\ref{3:1gx}) implies that
	\begin{eqnarray}
		H_\bullet(x)=\lambda \left\{\left(\a_{1}x+\b_{1}\overline{\a}_{1}\right)^{Q_{1}+Q_{2}+1}+\a_{2}\left(x+\b_{1}\right)^{Q_{1}+Q_{2}+1}\right\}.\end{eqnarray}
	for some  $\lambda\in\f^{\ast}$. Now using
	 $$f_\bullet(x)=x^{Q_{1}+Q_{2}+1}H_\bullet\left(x^{q-1}\right)=x^{Q_{1}+Q_{2}+1}H_\bullet\left(\overline{x}/x\right),$$
	we can obtain
	 \[f_\bullet(x)=\lambda\left[\left(\a_{1}\overline{x}+\b_{1}\overline{\a}_{1}x\right)^{Q_{1}+Q_{2}+1}+\a_{2}\left(\overline{x}+\b_{1}x\right)^{Q_{1}+Q_{2}+1}\right]. \]
	The above expression of $f_\bullet(x)$ can be further simplified. Writing $\b_{1}=\overline{\epsilon}/\epsilon$ for some $\epsilon\in\f^{\ast},$ and letting
	\begin{equation}
		u = \epsilon \a_{1}\overline{x}+\overline{\epsilon} \overline{\a}_{1}x,\quad
		v = \epsilon\overline{x}+\overline{\epsilon}x,			
	\end{equation}
	clearly, $u,v\in\h$ for any $x\in\f$, we can write $f(x)$ as
	 \[f(x)=\frac{\lambda}{\epsilon^{Q_{1}+Q_{2}+1}}[u^{Q_{1}+Q_{2}+1}+\a_{2}v^{Q_{1}+Q_{2}+1}]=L_{1}\circ P_1\circ L_{2}(x),\]
	where the maps $P_1,L_1,L_2$ are defined as
	\begin{eqnarray*}
		P_1(u,v) &=& \left(u^{Q_{1}+Q_{2}+1},v^{Q_{1}+Q_{2}+1}\right):\h^{2}\rightarrow \h^{2},\\
		 L_{1}(u,v)&=&\frac{\lambda}{\epsilon^{Q_{1}+Q_{2}+1}}\left(u+\a_{2}v\right):\h^{2}\rightarrow \f,\\
		 L_{2}(x)&=&\left(\epsilon\a_{1}\overline{x}+\overline{\epsilon}\overline{\a}_{1}x,\epsilon \overline{x}+\overline{\epsilon}x\right):\f\rightarrow\h^2.
	\end{eqnarray*} 	
	Clearly, $L_{1},L_{2}$ are linear functions. They are also permutations: $L_{1}$ is bijective since $\lambda\epsilon\ne 0$ and $\a_{2}\in\f\setminus\h$; $L_{2}$ is bijective since
	$$\det\begin{pmatrix}
		\epsilon\a_{1} &   \overline{\epsilon}\overline{ \a}_{1} \\
		\epsilon & \overline{\epsilon}\\
	\end{pmatrix}=\epsilon^{q+1}\left(\a_{1}+\overline{\a}_{1}\right)\ne 0.$$
	Thus $f_\bullet(x)$ is linearly equivalent to $P_1(u,v)$ over $\f.$
	
	\noindent {\bf Case 2: $m$ is even.} From 2 of Lemma \ref{3:clag}, there are degree-one rational functions $\rho(x), \sigma(x)$ over $\f$ defined from $\u$ to $\u$ such that $g_\bullet(x)=\rho^{-1}\circ x^{Q_{1}+Q_{2}+1-2r_\bullet}\circ \sigma$. By Lemma \ref{u2u}, we may write $\rho$ and $\sigma$ as
	\[\rho(x)=\frac{\overline{\a}_{2}x+\overline{\a}_1}{\a_{1}x+\a_{2}},\quad \sigma(x)= \frac{\overline{\b}_{2}x+\overline{\b}_{1}}{\b_{1}x+\b_{2}},\]
	where $\a_i,\b_i \in \f$ for $i=1,2$ and $\a_{1}\overline{\a}_{1} \ne \a_{2}\overline{\a}_{2}, \b_{1}\overline{\b}_{1} \ne \b_{2}\overline{\b}_{2}$. Moreover, using
	\begin{eqnarray*} \rho: \begin{array}{l}
			\frac{\overline{\a}_1}{\overline{\a}_2} \mapsto 0, \\
			\frac{{\a}_2}{\a_1} \mapsto \infty,
		\end{array} \quad \sigma: \begin{array}{l}
			\frac{\overline{\b}_1}{\overline{\b}_2} \mapsto 0, \\
			\frac{{\b}_2}{\b_1} \mapsto \infty,
		\end{array}
	\end{eqnarray*}  and noting the proof of {\bf Case 2: $m$ is even} of Lemma \ref{3:clag}, we have
	\begin{eqnarray} \label{3:abpro}
		\frac{\a_2}{\a_1},\frac{\b_2}{\b_1} \in \left\{\omega,\omega^2\right\} \subset \mathbb{F}_q,
	\end{eqnarray}
	where $\omega,\omega^2$ are the two distinct roots of $Q(x)=x^2+x+1$. Thus we have
	\[g_\bullet (x)=\frac{\a_{2}x+\overline{\a}_{1}}{\a_1x+\overline{\a}_{2}}\circ x^{Q_{1}+Q_{2}+1-2r_\bullet}\circ \frac{\overline{\b}_{2}x+\overline{\b}_{1}}{\b_{1}x+\b_{2}}.\]
	Expanding the above expression of $g_\bullet(x)$ we can obtain
	 \[g_\bullet(x)=\frac{\a_2\left(\overline{\b}_2x+\overline{\b}_1\right)^{Q_{1}+Q_{2}+1-2r_\bullet}+\overline{\a}_{1}\left(\b_1x+\b_{2}\right)^{Q_{1}+Q_{2}+1-2r_\bullet}}{\a_1\left(\overline{\b}_2x+\overline{\b}_1\right)^{Q_{1}+Q_{2}+1-2r_\bullet}+\overline{\a}_{2}\left(\b_1x+\b_{2}\right)^{Q_{1}+Q_{2}+1-2r_\bullet}}.\]
	Noting that in this case $g_\bullet(x)=\frac{N_\bullet(x)}{H_\bullet(x)}=\frac{\widetilde{N}_\bullet(x)}{\widetilde{H}_\bullet(x)}$ where $\gcd\left(\widetilde{N}_\bullet(x),\widetilde{H}_\bullet(x)\right)=1$, $H_\bullet(x)=Q(x)^{r_\bullet} \widetilde{H}_\bullet(x)$ and $N_\bullet(x)=Q(x)^{r_\bullet} \widetilde{N}_\bullet(x)$. Since $\deg \widetilde{N}_\bullet(x)=Q_1+Q_2+1-2r_\bullet$ for any $\bullet \in \{A,B,C\}$, Eq (\ref{3:1gx}) implies that
	\begin{eqnarray*}
		H_\bullet(x)=\lambda Q(x)^{r_\bullet} \left\{\a_1\left(\overline{\b}_2x+\overline{\b}_1\right)^{Q_{1}+Q_{2}+1-2r_\bullet}+\overline{\a}_{2}\left(\b_1x+\b_{2}\right)^{Q_{1}+Q_{2}+1-2r_\bullet}\right\} \end{eqnarray*}
	for some  $\lambda\in\f^{\ast}$. By (\ref{3:abpro}), we may write
	\[Q(x)=\left(x+\omega\right) \left(x+\omega^2\right)=\frac{1}{\b_1\overline{\b}_2}\left(\overline{\b}_2x+\overline{\b}_1\right)\left(\b_1x+\b_2\right).\]
	So we have
	\begin{eqnarray*}
		H_\bullet(x)=\frac{\lambda}{\left(\b_1\overline{\b}_2\right)^{r_\bullet}} \left\{\a_1\left(\b_1x+\b_{2}\right)^{r_\bullet}\left(\overline{\b}_2x+\overline{\b}_1\right)^{Q_{1}+Q_{2}+1-r_\bullet}+\overline{\a}_{2}\left(\overline{\b}_2x+\overline{\b}_1\right)^{r_\bullet}\left(\b_1x+\b_{2}\right)^{Q_{1}+Q_{2}+1-r_\bullet}\right\} \end{eqnarray*}
	Now using
	 $$f_\bullet(x)=x^{Q_{1}+Q_{2}+1}H_\bullet\left(x^{q-1}\right)=x^{Q_{1}+Q_{2}+1}H_\bullet\left(\overline{x}/x\right),$$
	we can obtain
	\begin{eqnarray*} f_\bullet(x)&=&c\cdot \left[\a_1\left(\b_1\overline{x}+\b_{2}x\right)^{r_\bullet}\left(\overline{\b}_{2}\overline{x}+\overline{\b}_{1}x\right)^{Q_{1}+Q_{2}+1-r_\bullet}+\overline{\a}_{2}\left(\overline{\b}_{2}\overline{x}+\overline{\b}_{1}x\right)^{r_\bullet}\left(\b_1\overline{x}+\b_{2}x\right)^{Q_{1}+Q_{2}+1-r_\bullet}\right] \\
		&=& L_1 \circ x^{Q_1+Q_2+1-r_\bullet} \overline{x}^{r_\bullet} \circ L_2,\end{eqnarray*}
	where $c:=\frac{\lambda}{\left(\b_1\overline{\b}_2\right)^{r_\bullet}} \in \f^*$ and $L_1,L_2$ are given by
	\[L_1(x)=c\left(\a_1\overline{x}+\overline{\a}_2x \right), \quad L_2(x)=\b_1\overline{x}+\b_2x. \]
	The maps $L_1,L_2: \f \to \f$ are linear permutations because $\a_1\overline{\a}_1\ne \a_2 \overline{\a}_2$ and $\b_1\overline{\b}_1 \ne \b_2\overline{\b}_2$. Hence $f_\bullet(x)$ is linear equivalent to $x^{Q_1+Q_2+1-r_\bullet} \overline{x}^{r_\bullet}$ as claimed by 2 of Lemma \ref{3:clag2}.
\end{proof}
Finally, by combining Lemmas \ref{3:gcd}--\ref{3:clag2}, we can now prove Theorems \ref{1:thmA}--\ref{1:thmC}.

\noindent \textbf{Proofs of Theorems \ref{1:thmA}--\ref{1:thmC}} For any $\bullet \in \{A,B,C\}$, we have $f_{\bullet}(x):=x^{Q_1+Q_2+1}H_{\bullet}\left(x^{q-1}\right)$. It is known that $f_\bullet(x)$ permutes $\f$ if and only if $\gcd\left(Q_1+Q_2+1,q-1\right)=1$, $H_\bullet(x)$ has no roots in $\mu_{q+1}$, and $g_\bullet(x)$ permutes $\mu_{q+1}$. The conditions such that $H_\bullet(x)$ has no roots in $\mu_{q+1}$ have been characterized in Lemma \ref{3:hxo}. By Lemma \ref{3:clag2}, if $\gcd\left(N_\bullet(x),H_\bullet(x)\right)=Q(x)^{r_\bullet}$ and $H_\bullet(x)$ has no roots in $\mu_{q+1}$, here the exact values of $r_\bullet$ were obtained in Lemma \ref{3:gcd}, if $m$ is odd, then $f_\bullet(x)$ is linear equivalent to $P_1(u,v)=\left(u^{Q_1+Q_2+1},v^{Q_1+Q_2+1}\right): \mathbb{F}_q^2 \to \mathbb{F}_q^2$, which is a permutation if and only if $\gcd\left(Q_1+Q_2+1,q-1\right)=1$. This proves 1 of Theorems \ref{1:thmA}--\ref{1:thmC}. On the other hand, if $m$ is even, then $f_\bullet(x)$ is linear equivalent to $P_2(x)=x^{Q_1+Q_2+1+m_\bullet(q-1)}: \f \to \f$, which is a permutation if and only if $\gcd\left(Q_1+Q_2+1+m\bullet(q-1), q^2-1\right)=1$. This proves 2 of Theorems \ref{1:thmA}--\ref{1:thmC}. Now the proof of Theorems \ref{1:thmA}--\ref{1:thmC} is complete.

	\section{Conclusion} \label{fin:cond}
	In this paper, we obtained three general classes of permutation pentanomials of the form $x^t+x^{r_1(q-1)+t}+x^{r_2(q-1)+t}+x^{r_3(q-1)+t}+x^{r_4(q-1)+t}$ over $\f$ where $q=2^m$ and their linear equivalence, which provided simple explanation as to why 14 families of the 17 permutation pentanomials discovered by Zhang et al. in \cite{zhang2024more} were permutations.  We raise two questions while writing this paper.
	\begin{itemize}
		\item[(1)] Can one find a simple explanation as to why the remaining 3 families among those 17 proved in \cite{zhang2024more} are permutations, and can one find extensions of these 3 families of permutation pentanomials over $\f$?
		
		\item[(2)] What about odd characteristic?
	\end{itemize}
	
	
	\bibliographystyle{amsplain}
	\bibliography{FX1}

\providecommand{\bysame}{\leavevmode\hbox to3em{\hrulefill}\thinspace}
\providecommand{\MR}{\relax\ifhmode\unskip\space\fi MR }
\providecommand{\MRhref}[2]{%
  \href{http://www.ams.org/mathscinet-getitem?mr=#1}{#2}
}
\providecommand{\href}[2]{#2}
\begin{thebibliography}{10}

\bibitem{ayad2014permutation}
M.~Ayad, K.~Belghaba, and O.~Kihel, \emph{On permutation binomials over finite
  fields}, Bull. Aust. Math. Soc. \textbf{89} (2014), 112--124.

\bibitem{BARTOLI2021101781}
D.~Bartoli and M.~Timpanella, \emph{A family of permutation trinomials over
  $\mathbb{F}_{q^{2}}$}, Finite Fields Appl. \textbf{70} (2021), 101781.

\bibitem{budaghyan2006new}
L.~Budaghyan, C.~Carlet, and A.~Pott, \emph{New classes of almost bent and
  almost perfect nonlinear polynomials}, IEEE Trans. Inf. Theory \textbf{52}
  (2006), 1141--1152.

\bibitem{carlet1998codes}
C.~Carlet, P.~Charpin, and V.~Zinoviev, \emph{Codes, bent functions and
  permutations suitable for des-like cryptosystems}, Des. Codes Cryptogr.
  \textbf{15} (1998), 125--156.

\bibitem{deng2019more}
H.~Deng and D.~Zheng, \emph{More classes of permutation trinomials with niho
  exponents}, Cryptogr. Commun. \textbf{11} (2019), 227--236.

\bibitem{ding2006family}
C.~Ding and J.~Yuan, \emph{A family of skew hadamard difference sets}, J. Comb.
  Theory, Ser. A \textbf{113} (2006), 1526--1535.

\bibitem{Ding2024}
Z.~Ding and M.~E. Zieve, \emph{On a class of permutation quadrinomials}, J.
  Algebra Appl. \textbf{0} (0), 2550075.

\bibitem{ding2023constructing}
\bysame, \emph{Constructing permutation polynomials using generalized redei
  functions}, arXiv preprint arXiv:2305.06322 (2023).

\bibitem{ding2023determination}
\bysame, \emph{Determination of a class of permutation quadrinomials}, Proc.
  Lond. Math. Soc. \textbf{127} (2023), 221--260.

\bibitem{gupta2024class}
R.~Gupta and A.~Rai, \emph{A class of permutation quadrinomials over finite
  fields}, Commun. Algebra \textbf{52} (2024), 1518--1524.

\bibitem{gupta2016some}
R.~Gupta and R.~K. Sharma, \emph{Some new classes of permutation trinomials
  over finite fields with even characteristic}, Finite Fields Appl. \textbf{41}
  (2016), 89--96.

\bibitem{Terry2024}
Chan~C. H. and Xiong M., \emph{Classification of a class of planar
  quadrinomials}, arXiv preprint arXiv:2404.14291 (2024).

\bibitem{hou2013class}
X.~Hou, \emph{A class of permutation trinomials over finite fields}, Acta
  Arith. \textbf{162} (2014), 51--64.

\bibitem{hou2014determination}
\bysame, \emph{Determination of a type of permutation trinomials over finite
  fields}, Acta Arith. \textbf{3} (2014), 253--278.

\bibitem{hou2015determination}
\bysame, \emph{Determination of a type of permutation trinomials over finite
  fields, {II}}, Finite Fields Appl. \textbf{35} (2015), 16--35.

\bibitem{hou2015survey}
\bysame, \emph{A survey of permutation binomials and trinomials over finite
  fields}, Contemp. Math. \textbf{632} (2015), 177--191.

\bibitem{laigle2007permutation}
Y.~Laigle-Chapuy, \emph{Permutation polynomials and applications to coding
  theory}, Finite Fields Appl. \textbf{13} (2007), 58--70.

\bibitem{li2017new}
K.~Li, L.~Qu, and X.~Chen, \emph{New classes of permutation binomials and
  permutation trinomials over finite fields}, Finite Fields Appl. \textbf{43}
  (2017), 69--85.

\bibitem{li2018new}
K.~Li, L.~Qu, and Q.~Wang, \emph{New constructions of permutation polynomials
  of the form $x^{r}h(x^{q-1})$ over $\mathbb{F}_{q^{2}}$}, Des. Codes and
  Cryptogr. \textbf{86} (2018), 2379--2405.

\bibitem{li2021several}
L.~Li, Q.~Wang, Y.~Xu, and X.~Zeng, \emph{Several classes of complete
  permutation polynomials with niho exponents}, Finite Fields Appl. \textbf{72}
  (2021), 101831.

\bibitem{LZ}
N.~Li and X.~Zeng, \emph{A survey on the applications of {N}iho exponents},
  Cryptogr. Commun. \textbf{11} (2019), no.~3, 509--548.

\bibitem{masuda2009permutation}
A.~Masuda and M.~Zieve, \emph{Permutation binomials over finite fields}, Trans.
  Amer. Math. Soc. \textbf{361} (2009), 4169--4180.

\bibitem{muller1986cryptanalysis}
W.~B M{\"u}ller and R.~N{\"o}bauer, \emph{Cryptanalysis of the dickson-scheme},
  Advances in Cryptology—EUROCRYPT’85: Proceedings of a Workshop on the
  Theory and Application of Cryptographic Techniques Linz, Austria, April 1985
  4, Springer, 1986, pp.~50--61.

\bibitem{ozbudak2023classification}
F.~{\"O}zbudak and B.~G. Tem{\"u}r, \emph{Classification of some quadrinomials
  over finite fields of odd characteristic}, Finite Fields Appl. \textbf{87}
  (2023), 102158.

\bibitem{rivest1978method}
R.~L. Rivest, A.~Shamir, and L.~Adleman, \emph{A method for obtaining digital
  signatures and public-key cryptosystems}, Commun. ACM \textbf{21} (1978),
  120--126.

\bibitem{shi2024ccz}
C.~Shi, J.~Peng, H.~Kan, and L.~Zheng, \emph{On ccz-equivalence between the
  bracken-tan-tan function and power functions}, Finite Fields Appl.
  \textbf{93} (2024), 102340.

\bibitem{sun2005interleavers}
J.~Sun and O.~Y. Takeshita, \emph{Interleavers for turbo codes using
  permutation polynomials over integer rings}, IEEE Trans. Inf. Theory
  \textbf{51} (2005), 101--119.

\bibitem{xu2018some}
G.~Xu, X.~Cao, and J.~Ping, \emph{Some permutation pentanomials over finite
  fields with even characteristic}, Finite Fields Appl. \textbf{49} (2018),
  212--226.

\bibitem{zhang2023more}
T.~Zhang, L.~Zheng, and X.~Hao, \emph{More classes of permutation hexanomials
  and pentanomials over finite fields with even characteristic}, Finite Fields
  Appl. \textbf{91} (2023), 102250.

\bibitem{zhang2024more}
T.~Zhang, L.~Zheng, and H.~Zhao, \emph{More classes of permutation pentanomials
  over finite fields with characteristic two}, Finite Fields Appl. \textbf{98}
  (2024), 102468.

\bibitem{zheng2022more}
L.~Zheng, B.~Liu, H.~Kan, J.~Peng, and D.~Tang, \emph{More classes of
  permutation quadrinomials from niho exponents in characteristic two}, Finite
  Fields Appl. \textbf{78} (2022), 101962.

\bibitem{zieve2009some}
M.~E. Zieve, \emph{{On some permutation polynomials over $\mathbb{F}_{q}$ of
  the form $x^{r}h(x^{(q-1)/d})$}}, Proc. Amer. Math. Soc. \textbf{137} (2009),
  2209--2216.

\bibitem{zieve2013permutation}
\bysame, \emph{Permutation polynomials on $\mathbb{F}_{q}$ induced from
  bijective redei functions on subgroups of the multiplicative group of
  $\mathbb{F}_{q}$}, arXiv preprint arXiv:1310.0776 (2013).

\end{thebibliography}
	
\end{document}